\documentclass[12pt]{article}
\usepackage{amsmath}
\usepackage{amssymb}
\usepackage{amsthm}
\usepackage{amscd}
\usepackage{amsfonts}
\usepackage{dsfont}
\usepackage{cases}
\usepackage[applemac]{inputenc}
\usepackage{enumerate}
\usepackage[all]{xy}
\usepackage{esint}
\usepackage{graphicx}
\usepackage{url}
\theoremstyle{plain}

\theoremstyle{definition}

\newcommand{\E}{\mathrm E}

\newcommand{\G}{\Gamma}

\newcommand{\vfi}{\varphi}

\newcommand{\s}{\sigma}

\begin{document}
\title{Numerical Aspects of Computing Possible Equilibria for Resource Dependent Branching Processes with Immigration}\date{}
\author{F. Thomas Bruss \\Universit\'e Libre de Bruxelles}
\maketitle
\begin{abstract} \noindent This article studies the stability of solutions of equilibrium equations arising in so-called resource dependent branching processes. We argue that these new models, building on the model already presented by Bruss (1984 a), refined and elaborated in Bruss and Duerinckx (2015) and now extended to allow immigration, are suitable to cope with specific properties of human populations. Our main interest is here to understand under which conditions immigration may lead to an equilibrium.
At the same time, we would like to advertize resource dependent branching processes as possibly the best models to study such questions.  

The equilibrium equations for the new models we obtain are clear and informative for several important stability questions. The goal of the study of the specific examples we provide is to see where the impact of immigration is most visible, and in how far increased efforts of integration can cope with dangers of instability. Moreover we discuss the advantages and a weaker point of our model, and also include a brief look at continuous-state, continuous-time
branching processes as an alternative.


\medskip
\noindent {\bf AMS 2010 Math. subj. classific.}: 60J85.

\medskip\noindent{\bf  Short Running title}: Immigration and Equilibria
\end{abstract}
\section{Motivation and related work}
The knowledge about different types branching processes and their use as a tools to study the growth of populations has grown immensely during the last five or six decades (see e.g. Haccou et al. (2015)). The motivation behind understanding specifically more about human populations and about the impact of immigration on the development of sub-populations, is that we see this understanding today as an urgent challenge.  

Shorter term problems dealing with human populations can often be studied with presently available data, and predictions are typically tackled through econometric models. On longer terms, however, this seems difficult. We can hardly predict the environment in which future generations will live, and how societies will adapt to them. Hence it is likely that more flexible models are needed, and Resource-dependent Branching Processes (RDBPs) on which we will focus, are models  which offer new features and flexibility.

\smallskip
First versions of such processes had already been introduced by the author at the S.P.A.-conference in 1982  (see Bruss (1984 a)) but  only  few papers on RDBPs have been published, and RDBPs are not very known. For the present paper we understand RDBPs as  they are defined in section 2 of Bruss and Duerinckx (2015), with an extension in Bruss (2016). These RDBPs lie in between the domain of controlled Branching processes (see the recent book by Gonz\'ales, del Puerto and Yanev (2018)) and branching processes with random environments,  as exemplified in Geiger et al. (2003), Hautphenne (2012)), and others.  

RDBPs can also be seen as extensions of so-called $\varphi$-branching processes (see Zubkov (1974)) since  in  these processes, in each generation, certain subsets of individuals are withdrawn  from reproduction, whereas those, who do reproduce, do so independently of each other.   $\varphi$-processes with random $\varphi$ were studied by Yanev (1976), and generalized in Bruss (1980). In RDBPs this control function $\varphi$ is a also a random function; moreover, in each new generation a different {\it policy} of the society can imply a different $\varphi,$  as explained in the next section. 

\smallskip
Although we stay in the present article focussed on RDBPs, we would like to draw here attention of the reader also to the broad concept of unification provided by Kersting (2017). Unlike our approach, Kersting's work does not aim specifically at studying human populations but offers different interesting approaches to convincing models for populations growing in a random environment.  

\section {Control in time}
RDBPs are local models, tailored to be part of a global model. They are seen
as elements of a sequence of probability models which are adapted in time. It is thus a {\it sequence of models} which is supposed to model the development of the population under consideration. At each updating time, the sequence becomes a specific RDBP, i.e. a well-defined probability model. 

The sequence itself is thus no probability model in a strict sense because the probabilistic prescription at time $t$ may only hold for one generation, and then change. It may sound confusing that such a vague definition should be sufficient to draw conclusions, but for certain {\it macro-economic} conclusions such as, for instance, the possibility of survival,  a precise formulation is not always necessary. 
Here is a simple example for this.

\smallskip
 Let $(\tilde Z_n)$ be a Galton-Watson type process with $\tilde Z_0=1$ with independent reproduction in which, given no extinction occurred before generation $k$, we have
	$$\tilde Z_{k+1}= D^{(k)}\big (\tilde Z_k), k=1, 2, \cdots,$$
where $D^{(k)}(n)$ denotes the random number of descendants generated by $n$ individuals which reproduce independently with mean $m_k.$ Then 
$$\{m_k>0 {\rm ~for~ all~ k}\}~ {\rm and}~\{m:=\liminf_k(m_k)>1\} \implies P(\lim_k \tilde Z_k=0)<1,$$ since $(\tilde Z_k)$ is, apart from at most finitely many exceptions, stochastically larger than a supercritical Galton-Watson process with reproduction mean $(m>1)$. 

Note that this implication is independent of the interpretation of the values $m_k$ as long as we maintain the hypothesis of independent reproduction. We can see the $m_k$ as averaged means of reproduction where "averaged" means averaged with respect to different classes of individuals, or averaged with respect to different constraints to which they are submitted.
As we shall see later, they decrease on average when society rules become more restrictive, or - and now we speak of immigration - they increase on average if new individuals with higher reproduction rates enter the home-population.
\subsection{Sequential update} The sequential update of local models within the sequence of models is supposed to be governed by the current rules of the society and the currently  observed pa-

\bigskip
~~{\bf Generation}: $...~\,t-2,.......~t-1,.........~t~~~\Big|\,......\,~t+1, ........~~t+2~...$

\medskip
~~{\bf Model}: $~~~~~~~...\,~\G^{[t-2]}~........~\G^{[t-1]}~.......~\G^{[t]}~\Big|~......~\tilde\G^{[t+1]}~~.......~\tilde \G^{[t+2]}~...$

\bigskip
\centerline {Figure 1}

 \noindent \begin{quote}In Figure 1, the process $\G^{[t]}$ stands for a RDBP defined by the parameters known at time $t$. The society examines whether the objective would be met if, for {\it all} $\tau >t,$ the processes $\tilde \G^{[\tau]}$ had the same probabilistic prescription as $\G^{[t]}.$
If not, the obligation principle forces the society to take  measures so that $\G^{[\tau]}$ would achieve this for $\tau \ge t+1.$\end{quote}
rameters. The control can only be specified through the objectives the society would like to meet,  among which the objective to {\it survive} is always defined to be the dominant one.  
The rules according to which the local models are adapted are supposed to follow the {\it science obligation principle.} This principle imposes that the probability the society tries to endow the population with parameters (i.e. creating an environment concerning consumption, production of resources, and reproduction), such that, if these parameters could be maintained forever, the population could survive forever.  (For more details we refer to section 7.1.1 in Bruss (2016).) 

The scheme in Figure 1 displays the controlled sequence, that is, the sequential update  of "local" probability models.
Thus the society's obligation principle leads at each time $t$ to a new local model for the next generation. We say here "leads to" instead of "implies" because there is no implication in the sense of a functional relationship
between successive models. The philosophy is that we cannot predict the future, and thus should not try to formulate any future guidelines for the society, and the idea is simply to maintain freedom of choice within the limits of those parameters not jeopardizing survival.

\subsection{Advantage of local models}
 The advantage of such a sequential  model is that the society's decision which  measures of control should be taken in the next generation can be derived in a well-defined probability model without imposing a probabilistic prescription for the whole future which is, as we know,  typically unknown. 
 
The idea is that the long-term objectives should govern the development of the population. 
All sub-populations studied in our forthcoming main model are supposed to follow the above scheme. However, we allow  that
the objectives may depend on all sub-populations at the same time.

\subsection{Content of this paper}
In Section 3 we introduce our assumptions and the terminology we need   to define the main model.

Section 4 deals with  numerical results for equilibrium equations. Although the solutions are nice in the sense that they are in general continuous as functions of the parameters, they can be  sensitive to small changes. One objective of the present addendum is to point out this sensitivity in a more detailed way. The reason is that, if our results are hoped to help decision makers to shape immigration policies, then they should see in which situations the risk of instability and possibly disastrous outcomes is high.

Section 5 discusses numerical methods. Several examples are given to show how sensitive equilibria
may be according to small changes of parameters or of the resource claim distribution functions.

In Section 6 we collect the main conclusions. This includes comments on strong parts and weaker parts of our model.
Section 7 finally provides an outlook to possible directions of further research.

We should add saying here, that the present article is still preliminary work. It will present ideas in the hope to raise interest, both in questions about certain impacts of immigration policies, as well as in Resource Dependent Branching Processes as suitable models to study such questions. The author thinks that this hope is justified, but he will not try here to go into all details.
\section {The main model}
\subsection{Terminology and assumptions}
The model consists of three discrete-time RDBPs, each
as defined in Bruss and Duerinckx (2015), section  2, denoting the effectives of three sub-populations. The members of these branching populations are called  {\it individuals}; these reproduce, produce resources, and consume, and in general with different parameters. 

Moreover, all individuals are submitted to the current  society rules, according to which they may, or may not, receive sufficient resources. If so, they are supposed to stay in the population; if not they will emigrate or protest in the sense that they  refuse reproduction.
Time is measured in terms of number of generations, where the $k$-th generation is defined as the interval 
$]k-1,k].$

The following assumptions are the same for all individuals within the {\it same} specified sub-population, i.e. in the same local model: 

\medskip
(i) Reproduction of individuals: This is a random variable  according to an identical law$\{p_k\}_{0,1,2, \cdots}$ where $p_k=P(\# {\rm ~descendants~ of~ an~ individual})$ with $0<p_0\le p_0+p_1<1$ and mean $0<m<\infty.$ Unless stated otherwise, reproduction of an individual is  independent of the reproduction of other individuals. (Since certain subsets of individuals may be excluded from reproduction, we have no Galton-Watson process (GWP).)
\medskip
\begin{figure}[ht]
\centering
 \includegraphics[width=0.6\textwidth, angle=0]{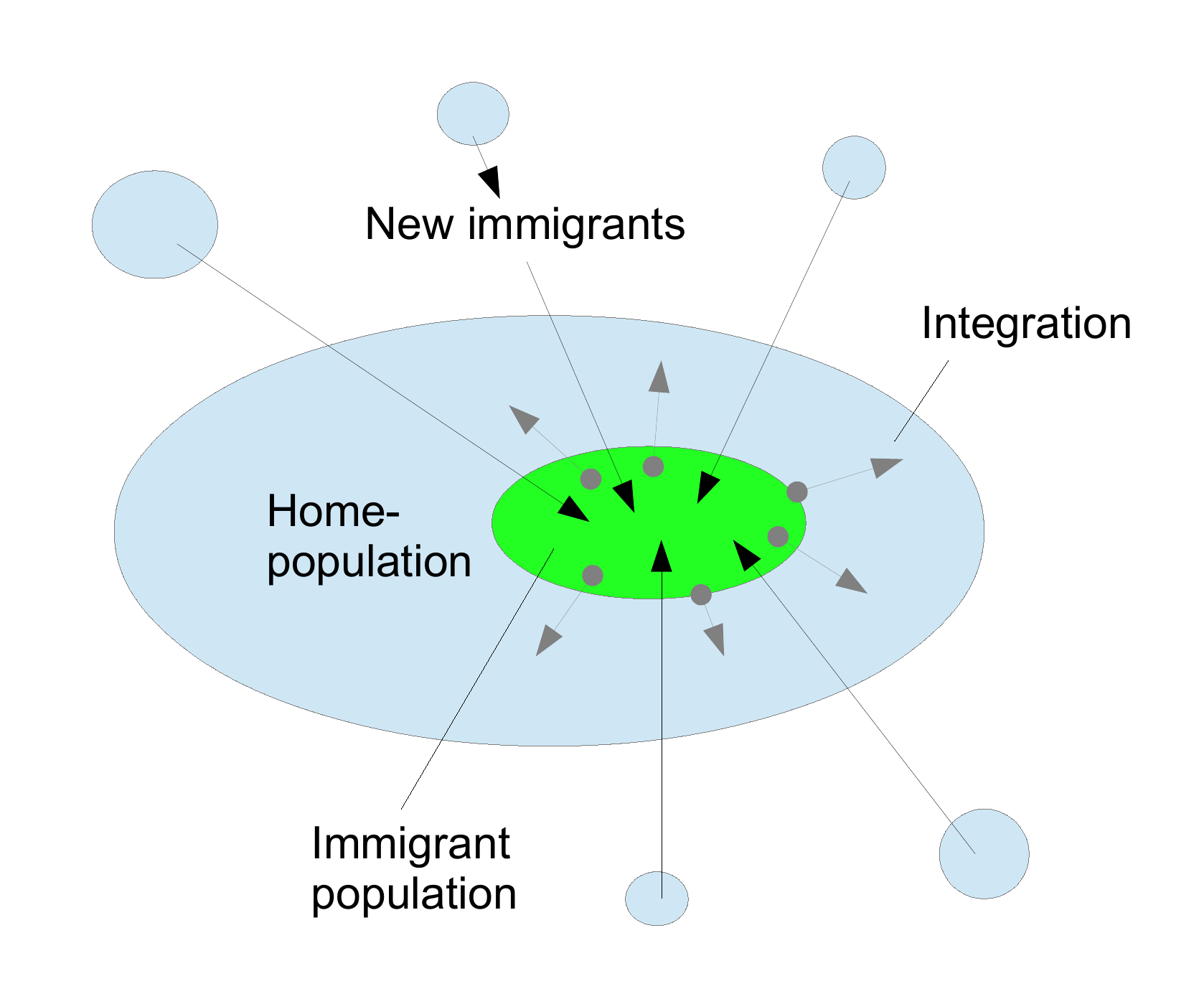}
 
 \centerline {Figure 2}
\begin{quote}Figure 2 shows the home-population and the immigrant population (smaller ellipse) living in cohabitation. In each generation a fraction of the latter integrates into the first one, and new immigrants arrive from other countries. \end{quote}
\end{figure}

(ii) Individual production of resources: It this paper, we only use the mean of this random variable $R$, denoted by $r.$

\medskip
(iii) Individual claims (consumption): Claims are supposed  to be random variables  governed by a continuous cumulative distribution function $F$, with finite second moments.

\medskip
(iv) Resources are supposed to be limited in each generation. The available resources are attributed to individuals according to currently valid rules. Individuals whose claims cannot be satisfied under these rules do not reproduce in the population.

\medskip\noindent
The set of current rules for the distribution of resources to individuals, to which we refer in (iv), is called the {\it current policy.} If the rules are supposed to stay the same over all generations we  speak of a {\it policy.}
For the precise definition of policies see Bruss and Duerinckx (2015), section 2. 

\subsection{The sub-populations}
In order to study the effect of immigration in a realistic way we have to consider at least three interacting sub-populations. We confine here to exactly three (see Figure 2). Correspondingly, our goal is to study a three-variate process $(\G_t)$ defined by
\begin{align} \G_t:=(\G_t^h, \G_t^i, I_t)_{~t=1,2, \cdots}\end{align} with the three components denoting the effectives of three sub-populations as explained below.

First the one we call
{\it home-population}, that is the one traditionally living in the hypothetically fixed country. Each individual belonging to that class is said to be in class $C_h,$ where $h$ is mnemonic for {\it home.} $C_h$ welcomes immigrants, and here we must make a distinction between, on the one hand, those who have settled already in the home-country but still maintain their identity and consider themselves as  being immigrants, and, on the other hand, new immigrants. 

The first class of immigrants is denoted by $C_i.$ In every generation, members of $C_i$ may change into members of $C_h,$ and in this case we say that they have {\it integrated} at time $t.$ The {\it integration rate} at time $t$ is denoted by $\varphi_t$, and we suppose  that the limit
$$ \lim_{t \to \infty} \varphi_t = \varphi$$ exists amost surely. Moreover, for the scope of the present article we can and do suppose that $\varphi_t\equiv \varphi$ 
To be definite we also suppose that integration happens at the beginning of the time interval $ ]t-1,t], t \in \{1, 2, \cdots\},$ and that their production, consumption and reproduction behavior of those individuals is at time $t$ exactly like members of the class $C_h.$ 

The second type of immigrants, i.e.\,the third sub-population, is called {\it new immigrants.} We denote this class by $C_{ni}.$ We assume that $C_{ni}$ is a transitory "class" in the sense that their offspring are by definition members of $C_i$ which then can integrate from the next generation onwards. It is no sub-population in the sense of the other two, as we shall see. We define 
\begin{align} I_t:= 
\#\{\rm{new~immigrants~arriving~at~time}~t\}\end{align}
The temporary distinction seems important for any realistic model because experience seems to confirm that  the difference in the needs and behavior of new immigrants is in general very different from that of  more established immigrants. 

\subsubsection{Terminology}
It is convenient to call Model I the model without immigrants, that is  $( \G_t):=(\G_t^h), {~t=1,2, \cdots}$. 

If $(\G_t)=(\G^h_t,\G_t^i)_{~t=1,2, \cdots}$, i.e. we have a model with two populations (one considered as home-population, the other as immigrant population) but no future new immigrants after some given time, then we speak of Model II. 

If moreover the stream of new immigrants may continue forever, we need the complete model (Modell III) defined in (1).

Accordingly, the model of Bruss and Duerinckx (2015) is Model I, thus just a special case of Model II,
which is itself a special case of Model III. In the Section 5 where we look at the numerical problems our graphs are mostly confined to Model II, since, except in special cases, the graphs for Model III would quickly become too complicated.

\subsection{Concentrating on the {\it wf}-policy}
Two special policies are singled out in Model I (Bruss and Duerinckx (2015)) for their importance. Only the so-called {\it weakest-first} policy ({\it wf}-policy), defined in subsection 3.2 of that paper
is essential in the present paper. The reason is that no society can survive unless it can survive under the {\it wf}-policy. 

The {\it wf}-policy stays equally important in the more complicated models, Model  II and the general model (Model III).

\smallskip

We first have to explain how the {\it wf}-policy is to be interpreted in the general model. Without survival there can be no equilibrium; thus survival of the sub-populations (home- and immigrant) is the ultimate necessary condition. (New immigrants are seen as a  temporary sub-population which will either integrate into the immigrant population, or else leave before). 

In our study, any rule to distribute resources is supposed to be strictly non-discriminating in the following sense: All claims of individuals are written down and then {\it merged} into one common list. Entries of this merged list are treated in the same way, independently of to which sub-population an individual (who submits the claim), belongs, and also independently of what the individual's direct ancestor consumed (claimed) or was able to contribute to the common total resource space. 

\smallskip\noindent
Hence the {\it wf}-policy can be summarized as follows:
\begin{quote} (i)~ In each generation, the society writes down the list of all individual claims for each sub-population, that is, if $n_\s$ is the length of the $\s$th list,
\begin{align}L^{(\s)}=\{X_1^{(\s)}, X_2^{(\s)}, \cdots, X_{n_\s}^{(\s)}\},\end{align}
where $X_j^{(\s)}$ denotes the $j$th claim from sub-population $\s.$ Here $\s\in\{1,2, \cdots, T\}$, where $T$ denotes the number of sub-populations. (In the present study, this mean $T=3.$) To simplify notation the generation number is not indicated.

\medskip
(ii)~The society orders than each list of claims $L^{(\s)}$ obtained in (i) into the corresponding list of increasing order statistics $$\tilde L^{(\s)}=\{X_{1, n_\s}^{(\s)}\le X_{2, n_\s}^{(\s)}\le  \cdots \le X_{n_\s, n_\s}^{(\s)}\},$$
and then merge all $T$ lists $\tilde L^{(\s)}$ into a single
list of increasing order statistics. This final list $L_Y,$
say, thus reads 
\begin{align}L_Y=\{Y_{1,z}\le Y_{2,z}\le \cdots \le Y_{z,z}\}, \end{align} where $ z=n_{\s_1}+n_{\s_2}+\cdots+n_{\s_T}.$

\medskip
(iii)~ Suppose that in the given generation $n,$ say, the current resource space is of size $s$, say. Let $N(z,s)$ be defined by
\begin{align}N(z,s):=\max\left\{0 \le k\le z: \sum_{j=1}^k
Y_{j,z}\le s\right\}.\end{align}\end{quote} Applying the {\it wf}-policy in generation $n$ means to  serve the $N(z,s)$ smallest claims
of the merged ordered list $L_Y$.

\medskip\noindent{\bf Remark}: Note that in (i) and (ii) we first
ordered the claims of the different sub-populations and then merged these ordered list into the final list
$L_Y.$ The outcome $L_Y$ would have been the same,  of course, if we had put all lists $L_\s$ together and then ordered all claims. Since all random variables within a given sub-population are i.i.d., we can however use more easily classical limit results to understand the asymptotic behavior of the  $\tilde L^{(\s)}.$  
\subsection{The sub-populations}
The  model we consider consists of three subpopulation. All three are now supposed to be RDBPs with their own characteristics. For ease of notation we use the mnemonic abbreviations $h,~i$ and $ni$ throughout for "home-population", "immigrant population" and "new immigrants. Thus for instance 
$F_h$ denotes the claim distribution function for individuals in $C_h, ~F_{ni}$ for individuals in $C_{ni},$ etc. Also, using $\delta$ as index notation, we recall that $r_\delta$ and $m_\delta$ denote for the corresponding classes the mean production of resources and the reproduction mean, respectively. 

\subsubsection{Cost functions}
The following integrals will play an important role in what will follow. They can be interpreted as the costs (in terms of spending resources) due to serving claims below an upper bound $t$ by the three different sub-populations.
\begin{align} \Psi_\delta(t):=\int_0^t x dF_\delta(x), ~\delta\in \{h,i, ni\}.\end{align} We define the row-vector function $\vec\Psi(x)$ by
$$\vec\Psi(x)=(\Psi_h(x), \Psi_i(x), \Psi_{ni}(x)).$$ 

 \medskip
 \begin{figure}[ht]
\centering
 \includegraphics[width=0.65\textwidth, angle=0]{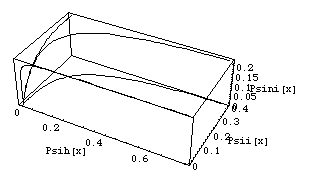}
 
 \centerline{Figure 3}
\begin{quote} Figure 3 shows for the  beta-distributions defined in (13) the graph of $\vec\Psi(x)$ running from bottom left (0,0,0) to top right. The labels are $\Psi_h(x)$, and correspondingly $\Psi_i(x)$ and $\Psi_{ni}(x).$ The graph shows also the projections of $\vec\Psi(x)$ on the three planes defined by one component set equal to zero \end{quote}
\end{figure}
These different cost functions  must be seen by keeping in mind that the corresponding mean production of members of these sub-populations, i.e. $r_h$, $r_i$ and $r_{ni}$ respectively, which are in general different. However, a priori, any sort of cost-benefit analysis within each sub-population is of little interest for any $t$ with $F_\delta(t)<1.$ because we suppose independence between individual claims and individual consumption. 

The $\Psi_\delta(t)$ defined by (6) are all continuously differentiable and strictly increasing $t$.

Suppose now that at generation $n$ the three sub-populations $C_h, C_i$ and $C_{ni}$ consisted of $g_h(n), g_i(n)$ and $g_{ni}(n)$ individuals, respectively. They all submit their claims for resources according to their respective claim distribution functions.
Suppose that the society  decides to serve if possible, (that is, if the available resource space is sufficient to do so) all claims
which do not exceed an upper bound $t$.
Then the society will face a total number of $g_h(n)+g_i(n)+g_{ni}(n)$ claims. Hence the total number of accepted claims $A(n,t),$ say, and the total  cost $K(n,t)$ (in terms of resources) under this policy satisfy, in expectation, $$\E(A(n,t))=\sum_{\delta \in \{h, i, ni\}}
g_\delta(n) F_\delta(t), ~~\E(K(n,t))=\sum_{\delta \in \{h, i, ni\}}
g_\delta(n) \Phi_\delta(t),$$ respectively. 
But then it is possible that the resources do not suffice to serve all claims not exceeding $t.$ 

Also, and in particular, we would like to understand what is the connection
between, on the one hand,  the upper bound $t$ for acceptable claims and, on the other hand, the total number of individuals which would be served under the {\it wf}-policy defined in Sub-section 3.3. Indeed, we recall that no population can survive unless it can survive as a {\it wf}-society.

\subsection {The BRS-inequality} 
 Remember that if we know the effectives of the sub-populations $g_h(n-1), g_i(n-1)$ and $g_{ni}(n-1)$ then we know also the expected total amount
of resources they leave for their descendants, namely
$$\sum_{\delta\in\{h,i,ni\}}g_\delta(n-1) \times r_\delta.$$
We now recall the following inequality (BRS-inequality), with which we can answer the last central question in the preceding Sub-section.

\bigskip \noindent {\bf Theorem} (Bruss and Robertson (1991), Steele (2016)) 
 \noindent Let $X_1, X_2, \cdots, X_n$ be positive random variables such that each $X_k$ has a absolute continuous distribution function $F_k,$ and let $N(n,s)$ be defined as in (5). Then
\begin{align} \E(N(n,s)) \le \sum_{k=1}^n F_k(\tau), \end{align} where  $\tau:=\tau(n,s)$ is a solution of\begin{align} \sum_{k=1}^n \int_0^\tau x dF_k(x)= s. 
\end{align} 
We can use this inequality to give an upper bound for
the total number of accepted claims. 

Recall that the society is supposed to be non-discriminating so that, in each generation, the policy applies to all three sub-populations. These have (in general) different claim distributions $F_h, F_i,$ and $F_{ni}$
and all their descendants submit their claims for resources. Thus society faces  lists of claims
which must all be treated according to the same rules.
If, in a given generation, $\nu_\delta$  denotes the number of claims submitted by the sub-population $\delta\in\{h,i,ni\}$ then the society must look at the ordered list of all claims as described in Subsection 3.3, that is 

$$
\begin{cases}
X_1^h, X_2^h,\, \cdots ,X_{\nu_h}^h& {\rm List ~1} \\
X_1^i, \,X_2^i,\, \cdots , X_{\nu_i}^i&{\rm List ~2}
\\
X_1^{ni}, X_2^{ni},\cdots, X_{\nu_{ni}}^{ni}&{\rm List ~3}
\end{cases} ~~\longrightarrow {\rm ~~merged~list.}
$$

\medskip
Then the Society orders the merged list in increasing order, yielding $$Y_1\le Y_2\le \cdots \le Y_{\nu_h+\nu_i+\nu_{ni}}.$$
Hence (7) and equality (8) in the BRS-inequality translate for the three classes of claims with distribution function  $F_\delta,  \delta \in\{h, i, ni\}$ as follows
 \begin{align}~~~~~~~n=\,&\nu_h+\nu_i+\nu_{ni},\\\E(N(n,s))\le &\,\nu_h F_h(\tau)+\nu_i F_I(\tau)+\nu_{ni}F_{ni}(\tau)
\end{align}
where $\tau:=\tau(n,s)$ is a solution of
\begin{align}\nu_h\int_o^\tau x dF_h(x)+\nu_i\int_o^\tau x dF_i(x)+\nu_{ni}\int_o^\tau x dF_{ni}(x)=s,\end{align} that is, according to (6),  a solution of the equation $$\nu_h\Psi_h(\tau)+\nu_i\Psi(\tau)+\nu_{ni}\Psi_{ni}(\tau) = s.$$

\section{Equilibria between sub-populations}
Since both the limiting integration rate and new immigrants in each generation have an influence on $(\G_t^h)$ and $(\G_t^i)$ the notation (1) should read more precisely
\begin{align}\left(\G^{(\vfi,I)}(t)\right)_{t=1, 2, \cdots}=\left(\Gamma_t^{h,(\vfi,I)}, \Gamma_t^{i,(\vfi,I)}, I_t\right)_{t=1, 2, \cdots}\end{align}    If $(I_t)\equiv 0$ and $\vfi$ defined after (1) the fraction of members of the 
immigrant-population which integrates into the 
home-population, then the process $\left(\G^{(\vfi,I)}(t)\right)$ is defined as the bi-variate process describing the effectives of $C_h$ and $C_i$ in (pure) co-habitation.

We are now ready to recall the definition of an equilibrium for the general model.

\bigskip
\noindent{\bf Definition 1} The stochastic process $(\G(t)^{(\vfi, I)})_{t=1,2, \cdots}$ defined in (4) is said to converge to an {\it equilibrium}, if there exists a value $\vfi\in [0,1]$ and a corresponding random variable $\alpha:=\alpha_\vfi$, such that conditioned on the survival of both $(\G_t^h)$ and  $(\G_t^i)$
\begin{align} \lim_{t \to \infty} \frac{\G_t^{i, \,(\vfi, I)}}{\G_t^{h, \,(\vfi, I)}}  = \alpha_{\vfi}~{\rm a.s.}\end{align} 
{\bf Remark 1} It is not hard to show that,  under realistic real-world assumptions for the claim distributions, the random variable $\alpha$ in Model II will allow for only few isolated possible values, depending on the initial numbers of effectives of the sub-populations. This hold also for Model III. Moreover, for this addendum it is justified to think of $\alpha$, if it exists,
as being constant.
 \subsection{Modeling interactions}
The interaction of $(I_t)$ with the two sub-processes can be modeled in many ways, and each model may have its own justification.
 We want to maintain more freedom in the model. 
 As we shall see below our approach is sufficiently flexible. 
 
 We can keep all factors different and will arrive at tractable equations whenever $I_t/\G^h_t$ conditioned on survival of the processes $(\G_t^h)$ and $(\G_t^i),$ can  be assumed  to satisfy 
\begin{align}\exists\, \ell_{ni}\ge 0\,:\, \frac{I_t}{\G_t^h} \to \ell_{ni} {\rm~almost~surely~as~}  t\to \infty. \end{align}
Note that the hypothesis of the existence of the limit $\ell_{ni}$ is not that  restrictive since $\ell_{ni}=0$ is permitted. The existence of $\ell_{ni}$ is also a realistic assumption in practice because at least one of the sub-populations will typically have some influence on how many new immigrants will be allowed to enter. Clearly, we do not have to specify which one, because  $$\ell_{\it ni}=\lim_{t\to \infty} \frac{I_t}{\G_t^h}~{\rm and} ~\alpha=\lim_{t\to \infty} \frac{\G_t^i}{\G_t^h} {\rm~ exists~ a.s.} $$
$$\implies \lim_{t\to \infty} \frac{I_t}{\G_t^i}=\frac{\ell_{ni}}{\alpha} ~a.s.$$
\subsection {The central equilibrium equation} 
 {\bf Theorem 1} ~Let $\vec\Psi(t)$ denote the row vector function $\left(\Psi_h(t), \Psi_i(t),\Psi_{ni}(t)\right)$ for $0\le t\le \infty$, with $\Psi_\delta(t), \delta \in \{h, i, ni\}$ as defined in (6), and let $\vec R$ denote the row vector of reproduction means $\left(r_h, \,r_i,\, r_{ni}\right)$. A necessary condition for the existence of an $\alpha$-equilibrium for the stochastic process 
 $(\G(t)^{(\vfi, I})_{t=1,2, \cdots}$ with $ni$-limit $\ell_{ni}$ is the existence of values $0< \vfi <1$, and $\tau>0$ solving the equation \begin{align}
\vec\Psi(\tau) \cdot
        \begin{pmatrix}
          m_h(1+\vfi\,\alpha) \\         
          m_i\Big((1-\vfi)\alpha+G_{ni}(\tau)\Big)\\\ell_{ni}
          \end{pmatrix}=\vec R \cdot  \begin{pmatrix}
          1+\vfi\,\alpha \\         
          (1-\vfi)\alpha  \\
        G_{ni}(\tau)
          \end{pmatrix}    
  \end{align}
where
\begin{align}G_{ni}(t):=F_{ni}(t) m_{ni}\ell_{ni}\left(F_h(t)m_h(1+\vfi \alpha)\right)^{-1},\end{align}
and where the claim distribution functions $ F_{h}, F_{i}$ and $F_{ni}$ and  the parameters\newline$m_i, m_h, m_{ni}, \ell_{ni} $ satisfy moreover the (double) constraint
  \begin{align}
F_h(\tau) m_h \Big(1+\vfi\alpha \Big)= F_i(\tau)\,\Big(m_i(1-\vfi)+G_{ni}(\tau)\Big)\ge \,1.                
  \end{align} 
  
  \noindent {\bf Remark 2.} Note the resource productivity means $r_h, r_i, r_{ni}$
  intervene implicitly in the constraints specified in (9) because $\tau, \varphi$ and $\alpha$ must satisfy equation (7).
  
 \subsection{Proofs and interpretations.}
 The proof of Theorem 1 is built on the Theorem of envelopment of Bruss and Duerinckx (2015), on the BRS-inequality and on classical limit results. To be brief, we here confine ourselves to give the essential steps. 

 \medskip\medskip
 {\bf Summarized proof of Theorem 1}:
 
 \smallskip
 (i) All three component processes in (1) are supposed to be RDBPs. It follows from the Theorem of envelopment (Bruss and Duerinckx (2015), page 341, Th. 4.13)
that no RDBP can survive unless it can survive under the {\it wf}-policy. Since survival of all sub-processes is a necessary condition
for the existence of an equilibrium, the starting point
is to think  in terms of the {\it wf}-policy. This will be used in (ii) and (iii).

 \smallskip
(ii) 
The merged list is a list of increasing order statistics of claims. Also in each generation all claims must be satisfied from resources taken in the same common resource space. With growing sub-populations its empirical distribution function  must converge almost surely. This explains why  there is a unique limiting threshold value $\tau$ in (11) for claims.  This is the value $\tau$ figuring as argument in all functions in equation (7). 

 \smallskip
(iii) The expected maximum total number
of claims which can be satisfied under the 
{\it wf}-policy is thus the maximal number of increasing order statistics of variables on the {\it merged} list of claims, the sum of which does not exceed the currently available resource space.  
An  upper bound for this expected number is obtained from the BRS-inequality (Steele (2016)), Bruss and Robertson (1991)).   

 \smallskip
(iv) Since reproduction
and resource production are i.i.d {\it within} each sub-population, we can (locally) apply the law of large numbers. This explains why the limiting equation for an equilibrium contains only the expected values of these, that is, $m_\delta$ respectively, $r_\delta,$ with $\delta\in\{h, i, ni\}.$ A Borel-Cantelli type argument together with Th.\,2.2 in Bruss and Robertson (1991,  pp 615-616), 
guarantees then that, conditioned on  the sub-populations going to infinity, the limiting inequality for the upper bound described in (iii) becomes a limiting equality. 

 \smallskip
(v) Using (iii) and (iv) we can see that the l.h.s.\,of (7) is the total expected consumption and the r.h.s. of (7)
the total availability of resources, however both {\it normalized per individual} in the home population. 
Normalization shows up through the strong law of large numbers by dividing
the BRS-equation by $\G_t^h$ and letting $\G_t^h$ tend to infinity.

 \smallskip
(vi) Finally, conditions 
(8) turn out to be the corresponding super-criticality conditions for the
processes $(\G_t^h)$ and $(\G_t^i).$ Since survival of both processes is a necessary condition, these conditions become compelling.

\smallskip
This completes the proof. \qed

\subsection{Implications}
It lies in the nature of the problem that the solution in Theorem 1 is complicated.  Although we tried to keep the central model as simply as possible the equation (7) with (8) and the constraints (9)
contains six parameters, three functions, and four variables, namely  $\tau, \varphi, \alpha$ and $ \ell_{ni}.$  All four variables are, in general,  implicit functions of each other.  In the next Section we discuss how to funnel our search for a numerical solution by looking at suitable projections.

\section{Numerical approach} Once the parameters $m_h, m_i,m_{ni}, r_h, r_i$ and $r_{ni}$, and the claim distribution functions $F_h, F_i$ and $F_{ni}$ are given, we can study the question whether an equilibrium, and at which level, will exist. We recall that if we have no independence within each su-population, then we can only hope for necessary conditions for survival of the subpopulations, and thus for the existence of an equilibrium. 

The main equilibrium equation with (7) and (8) and the constraints (9) depend  then on $\varphi, ~t, 
~\ell_{ni}$ and $\alpha.$ As already explained, we must fix  for a graphical presentation two of the four variables, and the question is with which we should begin.

The most efficient approach is probably to first select a main objective, and then to concentrate on this objective. It will likely circle around questions concerning feasible integration rates $\varphi$ and feasible values for an equilibrium.

For example, we may wonder 
whether an equilibrium at (more or less) a desired level is possible. Then we can solve (7) for $\alpha$ as a function of $\varphi, t$ and $\ell_{ni}.$ Below we see an example where this is done for an initial choice of $\ell_{ni},$ namely $\ell_{ni} = 0$. Then it becomes nice because (8) vanishes and the equation (7) becomes linear in $\alpha.$ 

However, our first interest may be quite different: the integration rate $\varphi$ may be constrained by an upper bound of real-world feasibility, and then the question of priority would be what range of possible values of $\alpha$ can be obtained by varying $\varphi$ in its feasible range. We will only show a few representative examples. 

Among others, one will show how sensitive the target values can be with rather small changes in the claim distributions and/or the  parameters.

\subsection{Examples}To keep the graphs as simple as possible, we choose all three claim distributions $F_h, F_i$ and $F_{ni}$ with aNumerical Aspects same support, $[0,1]$ say. For this it is convenient to choose three beta-distributions, as for instance
\begin{align}F_h :=\beta(6,2) , F_i:=\beta(2,3), F_{ni}:=\beta(2,7).\end{align}The other parameters we will use throughout in the examples (with the exception of Figure 8 where we exchange the parameters $m_h$ and $m_i$) are
\begin{align}m_h=3.5,~ m_i=2.8, ~r_h=2,~ r_i=0.6,\end{align}whereas we neither fix $\ell_{ni}$ nor $m_{ni}, r_{ni}$ in order to leave some room for discussing influence of these parameters.  The parameters chosen in (9) and (10) were {\it not} chosen to produce nice cases in which an equilibrium exists but actually more or less randomly, except that the graphs they produce should not hide the features we want to discuss.

\begin{figure}[ht]

\centering
 \includegraphics[width=0.55\textwidth, angle=0]{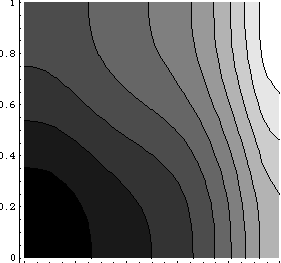}
 
 \centerline{Figure 4}
\begin{quote}Figure 4 shows the contour plot of 
$\Phi_{h,ni}(x,y)=\Psi_h(x)+\Psi_i(y)+\Psi_{ni}(x).$ The visible change from concavity to convexity of the contours for smaller $x$-values  $\varphi$-values for example would not be apparent in the corresponding plot of  $\Phi_{h,i}(x,y)=\Psi_h(x)+\Psi_i(x)+\Psi_{ni}(y).$\end{quote}
\end{figure}

Since the dot-product in the l.h.s. of (7) is in fact a weighted sum of the components of $\Psi(x),$  it can,
in view of understanding the sensitivity of possible solutions, be informative to look at either the sum $\Psi_h(x)+\Psi_i(x)+\Psi_{ni}(x)$, or more generally at the contour plot of functions depending on two threshold variables, as e.g. $\Psi_h(x)+\Psi_i(y)+\Psi_{ni}(x)$ or $\Psi_h(x)+\Psi_i(x)+\Psi_{ni}(y)$ with $(x,y)\in [0,1]^2.$ Here the effect of accumulated expected claim sizes below a {\it common} threshold become somewhat more visible
(by comparison), although it seems hard to give a clear interpretation of such partially accumulated sums.
\subsubsection {Visualizing relative costs}
Some information about a relative sensitivity of costs for different $t$-values can be gained by looking at certain contour plots. For example we may look at the contour plots
of the summarized cost by keeping the consumption level fixed for two sub-populations and vary the remaining one, as e.g.
$$\Phi_{h,i}(x,y):=\Psi_h(x)+\Psi_i(x)+\Psi_{ni}(y)$$
and the corresponding counterparts $\Phi_{h,ni}(x,y)$ and $\Phi_{i,ni}(x,y).$ For an example see Figure 3.

\begin{figure}[ht]

\smallskip

\centering
 \includegraphics[width=0.70\textwidth, angle=0]{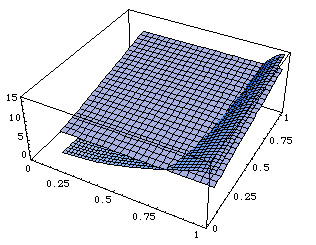}
 
 \centerline{Figure 5}
\begin{quote}Figure 5 shows the two surfaces describing the l.h.s. and the r.h.s. of equation (10) for the special case $\ell_{ni}=0.$ Here the "target" value $\alpha$ is put equal to $5.$ The r.h.s. is linear in both $\varphi$ and $t$ and thus the corresponding surface is the plane.
\end{quote}
\end{figure}

\subsection{Graphs of constraints for the equilibrium equation}
We now look at possibilities to graph the equation (7)
in an informative way. The possibly best start is to
graph it first for $\ell_{ni}=0$, that is for the case that
new immigrants become asymptotically negligible compared with the effectives in the home-population
and also compared with the immigrant population.
We then have also $G_{ni}(t) \equiv 0$ so that on both sides of equation (7) we have a dot-product of vectors of two remaining components. If, moreover, we fix the value of $\alpha,$ say, we obtain on each side a surface as a function of $t$ and $\varphi.$

The graph of these two surfaces as  a function of $t$ and $\varphi$
is shown in Figure 5 below for a fixed value $\alpha$. Since with $\ell_{ni}=0$ the r.h.s. of (7) is linear in $t$ and $\varphi,$ the r.h.s.-surface is thus the plane in this graph. 
\smallskip

All possible {\it candidates} for solutions $(t, \varphi)$ are the projection of the intersecting line between the plane and the curved surface on $[0,1]^2.$ Note that we cannot speak of solutions since the constraint qualifications given after equation (7) are not yet taken into account.

\subsection{Sensitivity of equilibria}

It is helpful to understand the sensitivity of a possible equilibrium $\alpha,$ in particular with respect to the limiting integration rate $\varphi.$ This is displayed in Figure 6.
\begin{figure}[ht]
\centering
 \includegraphics[width=0.70\textwidth, angle=0]{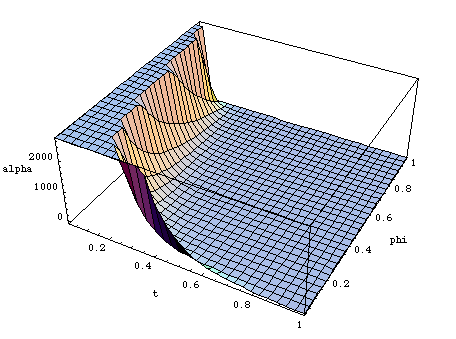}
 
 \centerline{Figure 6}
\begin{quote} Figure 6 exemplifies the sensitivity of possible equilibria. It shows for $\ell_{ni}=0$ the value of $\alpha$
as a function of $t$ and the limiting integration rate $\varphi.$  Note that the graph is truncated. It is the smaller possible values of $\alpha$ which are of interest (see (13)) if we want to maintain the natural meaning of an immigrant population integrates into the home-population.  
\end{quote}
\end{figure}
\begin{figure}[ht]
\centering{
 \includegraphics[width=0.75\textwidth, angle=0]{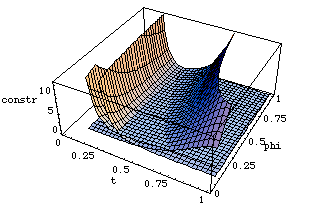}
 \centerline{Figure 7}}
 \end{figure}
\subsection{Interpretation of $\alpha$}

One must be careful in the interpretation of $\alpha.$
Only the values near the plane are of real-world interest. 

\medskip Moreover, Figure 7 above shows that it is often advisable to begin with the constraints for the existence of an equilibrium. We have chosen for Figure 7  the case of no new immigrants, implying in (15) in particular $\ell_{ni}=0.$ This allows us to display the fact that it is not necessarily the number of new immigrants per generation which might cause the sensitivity of solutions:

  For $\ell_{ni}=0$ we  have from 
 (16) and (17) that $G_{ni}(.)\equiv 0,$ so that (15)
 becomes the simpler equation
\begin{align}
\vec\Psi(\tau) \cdot
        \begin{pmatrix}
          m_h(1+\vfi\,\alpha) \\         
          m_i\Big((1-\vfi)\alpha\\
          \end{pmatrix}=\vec R \cdot  \begin{pmatrix}
          1+\vfi\,\alpha \\         
          (1-\vfi)\alpha  \\
          \end{pmatrix},    
  \end{align}
where $\vec\Psi(\tau)= (\Phi_h(\tau), \Phi_i(\tau))$ and $\vec R=(r_h, r_i),$
and where the claim distribution functions $ F_{h}, F_{i}$ and  the parameters $m_i, m_h, $ satisfy with equation (20) the double constraint
  \begin{align}
F_h(\tau) m_h \Big(1+\vfi\alpha \Big)= F_i(\tau)\,\Big(m_i(1-\vfi)\Big)\ge \,1.                
  \end{align} 
There are three surfaces shown in this graph, namely the l.h.s. of the constraint qualification (which is in Figure 7 the upper of the curved surfaces), the other curved surface for the corresponding r.h.s., and the plane $\pi[t,\varphi] \equiv 1.$ The two curved surfaces do intersect but their intersection lies {\it strictly} in the hidden part of the graph, that is, below the plane (critical for survival). Trying to have
have the home population and immigrant population steer towards an equilibrium would force them to intersect below the plane of criticality and thus to be wiped out both in the long run. Hence, no equilibrium can exist.

\noindent This contrasts the situation in Figure 8 shown below. Here we have exchanged $F_h$ and $F_i$ as well
as the parameters $r_h$ and $r_i$, and now the situation is different.

\begin{figure}[ht]

\centering
 \includegraphics[width=0.75\textwidth, angle=0]{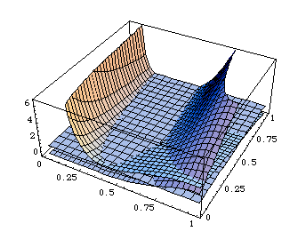}
 
 \centerline{Figure 8}
\end{figure}
\begin{quote} This figure (Fig 8) looks quite similar to Figure 7. Here $F_h$ and $F_i$ are interchanged. This implies the essential difference that, in the r.h.s. corner facing the viewer, the two surfaces do intersect, and that the curve of intersection lies strictly above the plane of criticality. The surface presenting $\alpha$ is not plotted in this graph. If it has points common
with this intersection curve, then each such point is a candidate
value for an equilibrium\end{quote}

\noindent
 Here we see that the two curved constraints surfaces have a curve of intersection above the plane of criticality.
  If the surface presenting $\alpha$ has points in common
with this curve, then each such point is a candidate
value for an equilibrium.

The valley between the two surfaces pointing to us describes the possible values  $\varphi$ and $\tau.$ An equilibrium should here de facto be  possible in reality because an integration rate $\varphi$  (r.h.s. axis) of about $25\%$ up to about $50\%$ per generation is not unrealistic. We do not see $\alpha$ in the same figure because the graph of $\alpha$ cannot be well combined in this graph without hiding important features.

 \section {Conclusions}
 
 The model which we presented  has several advantages. It is tractable and still  broad enough to deal with questions when sub-populations can reach an {\it equilibrium}, that is convergence of the ratio of respective effectives.  
The equations we derived are in terms of a policies to allocate resources, of mean natality and mean productivity rates, of the distribution functions of  {\it claims} of individuals of the home-population and the immigrant-population, respectively, as well as an {\it integration} parameter. 

There are cases where several candidates for an equilibrium may exist, but these are,
except in some cases of (unrealistic) perfect coincidences, very few. 

In real-world problems, where the number of individuals is large, of course,
it is typically clear from the current ratio of effectives, which candidate for an equilibrium will be the relevant one. This is why we suggested to think of $\alpha,$ if it exists, as being a constant.
\subsection{Importance of integration}
Our results show that conditions for the existence of an equilibrium prove to be severe, and also, how demanding the real world of immigration can be for politicians trying to make good decisions.  Secondly, we learn  {\it which} of the possible control measures stand out for combining feasibility and efficiency to reach an equilibrium, and how to recognize the steps one has to take towards controls.

 It is interesting to see that political support to help immigrants to integrate into the home-population, which is known to be important for a successful integration, is a very efficient control instrument for reaching an equilibrium. It is intuitive that it should be an efficient control because no other control has such a direct impact on the rate of change of the ratio of effectives. It is true that measures aiming to encourage changing the birth rates, say, would  be also quite effective, but this may be even harder to achieve than increasing the rate of integration. Birth rates are often deeply anchored in the sub-populations' respective traditions.
 
 However one must  interpret things correctly, in particular if the integration rate $\varphi$ becomes large. Suppose for example an equilibrium can be reached with an integration rate $\varphi=1/2,$ and this at level $\alpha=1/2.$ Thus in the limit (see (5)) we have one third of the limiting home-population being 
 immigrants with an expected integration time of $2$ generations. Suppose now that $\ell_{ni}>0.$ With ongoing immigration and turn-over rate $2$,   a randomly chosen individual  in the home-population is then, viewing his or her ancestors and  cultural background,  already much closer to a randomly chosen "immigrant" than to an "ancient" home-ancestor, and this long before the equilibrium is  reached. This is a fact and by no means an indication that the model may be unrealistic.

\subsection{Strong parts} Having pointed out the convincing parts of our model,
we can also argue that several simplifications in our model are inoffensive. 

\subsubsection{Accessible solutions for "realistic" hypotheses} In our model we can compute all necessary conditions for the existence of an equilibrium as "explicitly" as one can hope for, even if certain equations are bound to be implicit equations. From the numerically point of view there are no further difficulties to see and/or predict equilibria.

Moreover, the hypotheses we use can be well defended as being "realistic", at least in a preliminary sense. Much effort was indeed invested in trying to keep hypotheses defendable for human populations. We have hardly touched in this study the question of sufficient criteria for the existence of equilibria, but should add that, similarly as in the paper by Bruss and Duerinckx (2015) without immigration, several of the necessary conditions are also sufficient if the independence hypothesis for individual production and consumption can be maintained within each sub-population. 
All these properties constitute presumably the greatest advantage of our model.

\subsubsection{Asexual versus sexual reproduction}
Secondly, we should comment on why our model allows for an essential simplification, namely that we neglect the bi-sexual character of human populations. Indeed, bi-sexual reproduction can make in certain branching process models an essential difference compared to asexual reproduction, and this with respect to many typical questions. (For an overview about bi-sexual branching process models see e.g. Molina (2010)), and 
for a specific comparison see  Alsmeyer and R\"osler (2002)). 

However,  this makes no difference  for our model and our objectives because we can (luckily) confine our interest to asymptotic results conditioned on non-extinction, because, by definition, equilibrium questions are asymptotic features. In this case it is only the average unit reproduction mean per mating unit which counts, and not the mating law itself, as proved in
Bruss (1984 b); see (1), p. 916. This important simplification is also meaningful in several other branching models.

\subsubsection{Discrete time versus continuous time}
The discrete time setting, which is convenient for a clean definition of a merged list of claims,  should not make much of a difference either, since our findings depend only on the asymptotic behavior
of growth of the sub-populations. Without going into details here (see again Bruss and Duerincks (2015), the reason behind this is that, in order to reach an equilibrium,
the two important sub-populations must  behave asymptotically like Galton-Watson processes with reproduction mean $m_hF_h(\tau)>1$ and $m_iF_i(\tau)>1$ for some $\tau:=\tau(\alpha).$
Supercritical Galton-Watson processes have this property.  See for example Gonz\'ales et al. (2004, 2005). See also Hautphenne (2002) for another type of supercritical processes, namely decomposable branching processes.
\smallskip

\subsection{Weaker parts}

Now to a weaker part of our model, as conceived by the author. It concerns the implementation of our conclusions.

\smallskip
\subsubsection{Implementation in practice}
Viewing applicability, the sequential scheme shown in Section 2 may sometimes be hard to implement in practice. It is the {\it society's obligation principle} which is at stake.

Recall that this principle (see Subsection 2.1)  requires measures must be taken to assure that the society would survive forever if the updated model were valid forever. No lower bound is exacted for the survival probability. This fact itself is no problem since the survival probability can be kept arbitrarily close to one as seen from the {\it safe-haven}-property of the {\it wf}-policy in Bruss and Duerinckx (2015). However, the same paper shows also, that with the omni-present wish to increase the standard of living  in a RDBP the society may take survival risks by always searching the critical case. 

Since the parameters must be estimated in each generation for the next update, there is a danger that
parameters are frequently estimated too favorably for survival. But then the society obligation principle loses its strength as control mechanism. 
We have not spoken here about inference and/or estimation, but these tasks intervene at all stages.
So, for instance, Bruss and Slavtchova-Bojkova (1999) give a simple example showing that, in particular near criticality, inference must be treated with care when taking decisions. However, in our
model (where we try to take into account the human weakness of  aiming for a higher standard of living),
the near-to-critical case is likely to be the typical one.

\section{An outlook}
Trying to understand the development of human populations under immigration seems to be a major endeavor.
RDBPs are not yet sufficiently exploited, and one can probably go much further. However, RDBPs alone are not  likely to provide all answers for problems involving immigration. Hence it is good to see that there is much research going on in other new and independent directions.  We mention here the work of Barczy et al. (2018) on so-called aggregation of Galton-Watson-processes, and again the work of Kersting (2017) unifying the approach to random environments in branching processes.

Moreover, although most controlled Branching Processes are in discrete time,
there seems quite a potential for control in some continuous state/time counterparts of models, as exemplified by the paper of Li (2006) which gives  conditions for the weak convergence of  Galton-Watson branching processes with immigration to continuous-time, continuous-state branching processes with immigration.  Similarly, when hearing  at the Badajoz 2018-Conference on Branching Processes the interesting presentations of Bansaye et al. (2018), and Foucart (2018), the author, less familiar with the continuous state/time world of controlled BPs, wondered whether it is possible include in their settings  "resource dependence" and "society control" in a tractable form. 

In such instances of hearing about new approaches, the author feels  what quite a few mathematicians, who are also interested in applications, tend to say:
There is never too much theory for real applications.

\section*
{\bf Acknowledgement} The author would like to thank Guy Louchard (Universit\'e Libre de Bruxelles) for several helpful discussions concerning  aspects of numerical stability of equilibrium equations.

\section*{References}
 ~~~Alsmeyer G. and  Rösler U. (2002) {\it Asexual versus promiscuous bisexual Galton-Watson processes: The extinction probability ratio.}
 Ann. of Appl. Probab. Vol. 12, 125-142.

\smallskip
Bansaye V., Caballeroa M. E., M\'el\'eard S. (2018)
{\it Scaling limits for general finite dimensional
population models}, in Abstracts of the Badajoz Workshop on BPs (2018)

\smallskip
Barczy M.,  Ned\'enyi F. K. and  G. Pap (2018) {\it On aggregation of multitype Galton-Watson branching processes with immigration}, arXiv:1711.04099v2.

\smallskip
Bruss F. T. (1980) ~{\it A Counterpart of the Borel-Cantelli Lemma}, J. Appl. Prob., Vol. 17, 1094-1101.

\smallskip
Bruss F. T. (1984 a)~{\it Resource Depending Branching Processes} (Abstract), Stoch. Processes and Th. Applic.,  16 (1), 36-36.

\smallskip
Bruss F. T.  (1984 b) ~{\it A Note on Extinction Criteria for Bisexual Galton-Watson Processes}, J. Appl. Prob., Vol. 21, 915-919.

\smallskip
Bruss F. T. (2016) ~{\it The Theorem of Envelopment and Directives of Control in Resource
Dependent Branching Processes}, in Springer Lecture Notes in Statistics, Vol. 219 (I. M. del Puerto et al., Eds), 119-136.

\smallskip
Bruss F. T. and Robertson J. B. (1991) ~{\it 'Wald's Lemma' for Sums of Order Statistics of i.i.d. Random Variables},
Adv.  Appl. Prob., Vol 23, 612-623.

\smallskip
Bruss F. T. and Slavtchova-Bojkova M. (1999)~{\it On waiting times to populate an environment and
a question of statistical inference}, J. Appl. Prob., Vol. 36, 261-267.

\smallskip
Bruss F. T. and Duerinckx M. (2015) ~{\it Resource dependent branching processes and the envelope of societies},
 Ann. of Appl. Probab., Vol. 25, Nr 1, 324-372.
 
 \smallskip
 Geiger J., Kersting G. and Vatutin V.A. (2003)
 {\it Limit theorems for subcritical branching processes in random environment}, Ann. de l'Institut Henri Poincar\'e (B), Probability and Statistics,
Vol. 39, Issue 4, 593-620.

   \smallskip 
 Gonz\'alez M., Molina M., del Puerto I. (2004),
 {\it Recent results for supercritical branching processes with random control functions}, Pliska Studia Math. Bulgarica, Vol. 16, 43-54.
 
 \smallskip
   Gonz\'alez M., Martinez R., del Puerto I. (2005),
   {\it Estimation of the variance for a controlled branching process}, TEST, Vol. 14, 199-213.

\smallskip
Gonz\'ales M., Minuesa C., del Puerto I. (2016),
{\it Maximum likelihood estimation and Expectation-
Maximization algorithm for controlled branching processes}, Computational Statistics \& Data Analysis,
Vol. 93, 209-227. 

 \smallskip 
 Gonz\'alez M., del Puerto I.,  Yanev G. (2018) ~
 {\it Controlled Branching Processes}, Wiley, Volume 2.

\smallskip
Haccou P.,  Jagers P. and  Vatutin V.A. (2005) ~{\it Branching Processes - Variation, Growth and Extinction of Populations}, Cambridge University Press.

\smallskip
Hautphenne S. (2012)
{\it Extinction probabilities of supercritical decomposable branching processes},
 J.  Appl. Prob.,  Vol 49, 639-651.

\smallskip
Kersting G. (2017), {~\it A unifying approach to branching processes
in varying environment},  ~arXiv:1703.01960v6.

\smallskip

Li  Z. (2006),  {\it A Limit Theorem for Discrete Galton-Watson Branching Processes with Immigration}, J. Appl. Prob., Vol. 43, Issue 1, 289-295.

\smallskip
Molina M. (2010) ~{\it Two-sex branching process literature}, in Lecture Notes in Statistics-Proc. 196, Springer-Verlag 
(Eds. M. Gonzales, I.M. del Puerto, R. Martinez, M. Molina, M. Mota, A. Ramos), 279-291.

\smallskip
Steele J. M. (2016) ~{\it The Bruss-Robertson Inequality: Elaborations, Extensions, and Applications}, Math. Applicanda, Vol. 44,  No 1, 3-16.

\smallskip
Yanev N. M. (1976)~{\it Conditions for degeneracy of $\phi$-branching processes with random $\varphi$}, Theory Prob. Applic.,  Vol. 20, 421-427.

\smallskip
Zubkov A.M. (1974), {Analogies between Galton-Watson processes and $\varphi$-branching processes}, Theory Prob. Applic., Vol. 19, 319-339.

\bigskip\bigskip
\noindent Author's address:

\medskip
\noindent
F. Thomas Bruss\\
~~~~~Universit\'e Libre de Bruxelles\\
~~~~~D\'epartement de Math\'ematique\\
~~~~~Campus Plaine, CP 210\\
~~~~~B-1050 Brussels, Belgium

\bigskip
\noindent Phone:  ++32-2-650-5893;  

~~~~~~++32-472-448801.

\end{document}